\documentclass[12pt,a4paper]{article}

\usepackage{latexsym,amsfonts,amsmath,amsthm,graphicx,amssymb,xcolor}

\usepackage[colorlinks=true,linkcolor=blue,citecolor=blue]{hyperref}

\newcommand{\satop}[2]{\stackrel{\scriptstyle{#1}}{\scriptstyle{#2}}}

\def\bs0{\boldsymbol{0}}

\allowdisplaybreaks

\begin{document}

\title{A fast Fourier transform method for computing the weight enumerator polynomial and trigonometric degree of lattice rules}

\author{Josef Dick\footnote{School of Mathematics and Statistics, UNSW, Sydney, Australia; email: josef.dick@unsw.edu.au;  The author is supported by a Queen Elizabeth 2 Fellowship of the Australian Research Council.}}

\date{\today}
\maketitle

The weight enumerator polynomial $W_{(g_1,\ldots, g_s)}$, which has previously been studied in coding theory~\cite{DS02, Si01, Sid, vL92}, association schemes~\cite{MS99}, orthogonal arrays~\cite{T08}, spherical $t$-designs~\cite{Bannai} and  digital nets \cite{DM12, DS02, T08}, of a lattice rule~\cite{niesiam, SJ}
\begin{equation*}
\frac{1}{N} \sum_{n=0}^{N-1} f\left(\left\{\frac{n (g_1,\ldots, g_s)}{N}\right\}\right),
\end{equation*}
where $\{x\}$ denotes the fractional part of a nonnegative real $x$, $N > 1$ is an integer, and $g_i \in \{1, 2, \ldots, N-1\}$, defined by
\begin{equation*}
W_{(g_1,\ldots, g_s)}(z) := \sum_{a=0}^{ds} z^a M_{(g_1,\ldots, g_s)}(a),
\end{equation*}
with
\begin{align*}
M_{(g_1,\ldots, g_s)}(a) := \#\{ & (k_1,\ldots, k_s) \in \{-d,-d+1,\ldots, d\}^s: \\ & k_1g_1+\cdots + k_sg_s \equiv 0 \pmod{N}, |k_1|+\cdots + |k_s| = a \},
\end{align*}
can be written as
\begin{align*}
W_{(g_1,\ldots, g_s)}(z) = & \sum_{\satop{(k_1,\ldots, k_s) \in \{-d,\ldots, d\}^s}{k_1g_1+\cdots + k_sg_s \equiv 0 \pmod{N}}} z^{|k_1|+\cdots + |k_s|} \\ = & \sum_{(k_1,\ldots, k_s) \in \{-d,\ldots, d\}^s} z^{|k_1|+ \cdots + |k_s|} \frac{1}{N} \sum_{n=0}^{N-1} \mathrm{e}^{2\pi \mathrm{i} n (k_1g_1+\cdots + k_sg_s)/N} \\ = &
\frac{1}{N}\sum_{n=0}^{N-1} \prod_{j=1}^s \sum_{k_j=-d}^d z^{|k_j|} \mathrm{e}^{2\pi \mathrm{i} n k_j g_j/N},
\end{align*}
where the last expression can be used to compute $W_{(g_1,\ldots, g_s)}$ in $\mathcal{O}(Nds^2)$ operations, which can be reduced to $\mathcal{O}(Nds)$ operations if one wants to compute the trigonometric degree $\rho$ \cite{CL, LS} of a lattice rule
via the equality
\begin{equation*}
\rho = \left\{\begin{array}{rl} 0 & \mbox{if } M_1 \neq 0, \\  \max\{a: 1 \le a \le d, M_1=\cdots = M_a = 0\} & \mbox{otherwise}, \end{array} \right.
\end{equation*}
where the trigonometric degree is defined by
\begin{equation*}
\frac{1}{N} \sum_{n=0}^{N-1} \mathrm{e}^{2\pi \mathrm{i} n (k_1g_1 + \cdots + k_sg_s)/N} = \int_0^1 \cdots \int_0^1 \mathrm{e}^{2\pi\mathrm{i} (k_1x_1+\cdots + k_sx_s)} \,\mathrm{d}x_1\cdots \,\mathrm{d} x_s = 0,
\end{equation*}
for all $(k_1,\ldots, k_s)\in \mathbb{Z}^s \setminus \{(0,\ldots, 0)\}$ with $|k_1|+\cdots + |k_s| \le \rho$, since this is equivalent to the maximum over all $\rho$ such that
\begin{equation*}
k_1g_1+\cdots + k_sg_s \not \equiv 0 \pmod{N} \quad \mbox{for all } (k_1,\ldots, k_s) \in \mathbb{Z}\setminus\{(0,\ldots, 0)\} \mbox{ and } |k_1|+\cdots + |k_s| \le \rho.
\end{equation*}


\begin{thebibliography}{10}

\bibitem{Bannai} E. Bannai, On the weight distribution of spherical t-designs. European J. Combin., 1, 19--26, 1980.

\bibitem{CL} R. Cools and J. N. Lyness, Three- and four-dimensional K-optimal lattice rules of moderate trigonometric degree. Math. Comp., 70, 1549--1567, 2001.

\bibitem{DM12} J. Dick and M. Matsumoto, On the fast computation of the weight enumerator
polynomial and the $t$ value of digital nets over finite abelian groups. In preparation.

\bibitem{DS02} S. T. Dougherty and M. M. Skriganov, MacWilliams duality and the Rosenbloom-Tsfasman metric. Mosc. Math. J., 2(1), 81--97, 199, 2002.

\bibitem{LS} J. N. Lyness and T. S{\o}revik, Five-dimensional K-optimal lattice rules. Math. Comp., 75, 1467--1480, 2006.

\bibitem{MS99} W. J. Martin and D. R. Stinson, Association schemes for ordered orthogonal arrays and $(T,M,S)$-nets. Can. J. Math., 51, 326--346, 1999.

\bibitem{niesiam} H. Niederreiter, {\it Random number generation and quasi-Monte Carlo methods.} CBMS-NSF Regional Conference Series in Applied Mathematics, 63. Society for Industrial and Applied Mathematics (SIAM), Philadelphia, PA, 1992.
    
\bibitem{Si01} I. Siap, The complete weight enumerator for codes over $\mathcal{M}_{n,s}(\mathbb{F}_q)$. In: Cryptography and Coding, Lecture Notes in Computer Science, vol. 2260, pp. 20--26. Springer, Berlin (2001).

\bibitem{Sid} V. M. Sidel'nikov, The spectrum of weights of binary Bose-Chaudhuri-Hocquenghem codes. Problemy Pereda\v{c}i Informacii, 7, 14--22, 1971.

\bibitem{SJ} I. H. Sloan and S. Joe, Lattice methods for multiple integration. Oxford University Press, Oxford, 1994.

\bibitem{T08} H. Trinker, A simple derivation of the MacWilliams identity for linear ordered codes and orthogonal arrays. Des. Codes Cryptogr., 50, 229--234, 2009.

\bibitem{vL92} J. H. van Lint, Introduction to Coding Theory. Graduate Texts in Mathematics, vol. 86, 2nd edn. Springer-Verlag, Berlin (1992).


\end{thebibliography}
\end{document}